\documentclass{gtart}
\usepackage{amsmath, amssymb, %graphicx, %color, %epsfig, 
verbatim, pinlabel
}

\usepackage{amscd} \usepackage{amssymb} 

\newtheorem{thm}{Theorem}[section]

\newtheorem{cor}[thm]{Corollary}
\newtheorem{lem}[thm]{Lemma}
\newtheorem{prop}[thm]{Proposition}

\newtheorem{conj}[thm]{Conjecture}

\newtheorem{rem}[thm]{Remark}

\numberwithin{equation}{section}

\newcommand{\bfz}{{\mathbb {Z}}}

\newcommand{\oma}{\mathbf a}
\newcommand{\omb}{\mathbf b}

\newcommand{\omz}{\mathbf z}

\newcommand{\Z}{\mathbb Z}

\newcommand{\R}{\mathbb R}
\newcommand{\bfr}{\mathbb R}
\newcommand{\bfc}{\mathbb C}

\newcommand{\ra}{\rightarrow}

\DeclareMathOperator{\interior}{int}

\DeclareMathOperator{\genus}{genus}

\begin{document}

\title{Symplectic surgeries and normal surface singularities}

\author{David T. Gay}
\address{Department of Mathematics and Applied Mathematics\\ 
University of Cape Town\\
Private Bag X3, Rondebosch 7701\\
South Africa}

\author{Andr\'{a}s I. Stipsicz}
\address{R{\'e}nyi Institute of Mathematics\\
Re{\'a}ltanoda utca 13--15, Budapest, Hungary and \\
Mathematics Department, Columbia University\\
2990 Broadway, New York, NY 10027 }

\email{David.Gay@uct.ac.za, stipsicz@math-inst.hu}

\begin{abstract} 
We show that every negative definite configuration of symplectic surfaces in a symplectic $4$--manifold has a strongly symplectically convex neighborhood. We use this to show that, if a negative definite configuration satisfies an additional negativity condition at each surface in the configuration, and if the complex singularity with resolution diffeomorphic to a neighborhood of the configuration has a smoothing, then the configuration can be symplectically replaced by the smoothing of the singularity. This generalizes the symplectic rational blowdown procedure used in recent constructions of small exotic $4$--manifolds.
\end{abstract}
\maketitle

\section{Introduction}
\label{sec:first}
Most of the recent examples in smooth 4--manifold topology have been constructed
using the following ``cut-and-paste'' scheme: Suppose that the 
smooth closed 4--manifold $X$ is decomposed along the embedded 3--manifold $Y$
as 
\[ 
X=X_1\cup _Y X_2
\]
where $X_1, X_2$ are codimension--0 submanifolds of $X$ with $\partial
X_1 = -\partial X_2 = Y \neq \emptyset$. Suppose furthermore that
$Z_1$ is a smooth 4--manifold with boundary $\partial Z_1$
diffeomorphic to $Y = \partial X_1$.  Then a new 4--manifold
\[
Z=Z_1\cup _Y X_2
\]
can be constructed by cutting $\interior(X_1)$ out of $X$ and gluing
$Z_1$ back in. The topological type of $Z$ might also depend on the
gluing diffeomorphism $\varphi \colon \partial Z_1\to Y$, but for
simplicity we will suppress this dependence in the notation. For
example, if $X_1$ is a tubular neighborhood of a torus of
self--intersection 0 and $Z_1=D^2\times T^2$ then appropriate choices
of $\varphi$ give (generalized) logarithmic transformations and
Luttinger surgeries.

The most important topological invariants of a closed smooth
4--manifold are the fundamental group $\pi _1$, the Euler
characteristic $\chi $ and the signature $\sigma$.  In fact, in the
simply connected case $\chi $ and $\sigma$ essentially determine the
smooth 4--manifold up to homeomorphism \cite{Fr}.  The change of $\chi
$ and $\sigma$ can be very easily determined in a cut-and-paste
operation, since these quantitites are additive, while the fundamental
group can be computed using the Seifert--Van Kampen theorem. The
determination of the smooth structure is, however, much more
complicated. The most sensitive smooth invariant, the Seiberg--Witten
function
\[
SW _Z \colon H^2 (Z; \bfz ) \to \bfz 
\]
is very hard to compute in general, and although a TQFT-type theory
(the monopole Floer homology \cite{KM}) has been developed to compute
the Seiberg--Witten invariants of the result of a cut-and-paste
construction, such computations are extremally challenging in
practice. Partial knowledge of $SW_Z$ is provided by Taubes' famous
theorem \cite{Tau}, stating that $SW_Z (c_1(Z, \omega ))$ is $\pm 1$
provided $\omega \in \Omega ^2 (Z)$ is a symplectic form on $Z$ (and
$b_2^+(Z)>1$). Therefore we are particularly interested in
cut-and-paste constructions which can be performed within the
symplectic category.

In this paper we will consider the following special case of the above
cut-and-paste construction: Suppose that $C=C_1 \cup \ldots \cup C_n
\subset (X, \omega )$ is a collection of closed symplectic
2--dimensional submanifolds of the closed symplectic 4--manifold $(X,
\omega )$, intersecting each other $\omega$--orthogonally according to
the plumbing graph $\Gamma$.  Recall that each vertex $v$ of the
plumbing graph $\Gamma$ corresponds to a surface, hence is decorated
by two integers, the genus $g_v$ and the homological square (or
self--intersection) $s_v$ of the surface, and two verices are
connected by $n\geq 0$ edges if and only if the corresponding surfaces
intersect each other transversely in $n$ (positive) points. We will
denote the number of edges emanating from a vertex $v$ by $d_v$.  Let
$X_1$ be a tubular neighborhood $\nu C$ of the configuration $C=C_1
\cup \ldots \cup C_n$. Assume that $\Gamma $ is negative definite
(i.e. the corresponding intersection form is negative definite), and
consider a normal surface singularity $(S_{\Gamma }, 0)$ with
resolution graph $\Gamma$.  (It is a result of algebraic geometry
\cite{grau} that such $(S_{\Gamma }, 0)$ exists for every negative
definite $\Gamma $, although the analytic structure on $(S_{\Gamma},
0)$ might not be uniquely determined by $\Gamma$.)  Suppose finally
that $Z_1$ is the Milnor fiber of a smoothing of the singularity
$(S_{\Gamma }, 0)$.  Depending on $(S_{\Gamma}, 0)$, such smoothing
may or may not exist. For example, if $(S_{\Gamma}, 0)$ is a
hypersurface singularity (given by a single equation), or more
generally it is a complete intersection (cf. Section~\ref{sec:gen}),
then such smoothing always exists.  The main result of this paper is:

\begin{thm}\label{t:main}
Suppose that $\Gamma$ is a negative definite plumbing graph which either
\begin{enumerate}
\item is a tree and has $g_v=0$, $-s_v-d_v\geq 0$  for all vertices; or
 \item $-s_v> d_v +2g_v$ holds for every vertex $v$.
\end{enumerate}
Suppose furthermore that $C=C_1 \cup \ldots \cup C_n \subset (X,
\omega )$ is a collection of closed symplectic 2--dimensional
submanifolds of the closed symplectic 4--manifold $(X, \omega )$,
intersecting each other $\omega$--orthogonally according to the
plumbing graph $\Gamma$. Let $(S_{\Gamma}, 0)$ denote a singularity with
resolution graph $\Gamma$ and $Z_1$ the Milnor fiber of a smoothing of
$(S_{\Gamma}, 0)$. If $X_1\subset X$ is a closed tubular neighborhood of 
$C$ in $X$, then the 4--manifold 
\[
Z=Z_1\cup _Y (X -
\interior(X_1))
\]
(with a suitable, naturally chosen gluing diffeomorphism $\varphi$
specified later) admits a symplectic structure $\omega _Z$, which can
be assumed to agree with the given symplectic structure $\omega $ on
$X - \interior(X_1)$.
\end{thm}

One way of interpreting this result is the following: Consider the
singular 4--manifold $X^{sing}$ we get by collapsing $C$ to a
point. If the singularity of $X^{sing}$ is diffeomorphic to a
holomorphic model admitting a smoothing, and $\Gamma$ satisfies one of
the additional hypotheses given in the theorem, then this smoothing
can always be ``globalized'' in the symplectic category.  Notice that
we do not require the singular point to have a holomorphic model in
$X^{sing}$ as in \cite{McW2} (where the analytic structure near the
singular point is also assumed to be modeled by the holomorphic
situation) --- we just require the existence of a diffeomorphism. For
``globalizing'' local deformations in the holomorphic category in a
similar context, see \cite{LPark}.

According to \cite{CNP}, the link $Y=\partial Z_1$ of the singularity
$(S_{\Gamma }, 0)$ given by the (negative definite) plumbing graph
$\Gamma$ admits a unique (up to contactomorphism) Milnor fillable
contact structure $\xi_M$, for which $Z_1$ (with its Stein structure
originating from the deformation) provides a Stein filling. In fact,
our proof will not use the fact that $Z_1$ is a smoothing of
$(S_{\Gamma}, 0)$. Instead, we will rely on the fact that $Z_1$ admits
a symplectic structure $\Omega$ such that $(Z_1,\Omega)$ is a strong
symplectic filling of $(Y,\xi_{M})$. For this reason the chosen
analytic structure on $(S_{\Gamma }, 0)$ is not relevant.

For the convenience of the reader, below we summarize the strategy we
will use in the proof of Theorem~\ref{t:main}.  First we will show
that the union $C \subset (X,\omega)$ of the symplectic surfaces (of
arbitrary genera, intersecting each other $\omega$--orthogonally and
according to the negative definite graph $\Gamma$) in the symplectic
4--manifold $(X, \omega )$ admits a compact $\omega$--convex
neighborhood $U_C$.  This will be achieved by producing a model
symplectic $4$-manifold $(X_\Gamma,\omega_\Gamma)$ containing a
configuration $C_\Gamma$ of symplectic surfaces (intersecting each
other $\omega_\Gamma$--orthogonally and according to $\Gamma$), with the
same areas and genera as the surfaces in $C$ and with a neighborhood
system of $\omega_\Gamma$--convex neighborhoods of $C_\Gamma$, such
that any neighborhood $\nu C_{\Gamma }$ of $C_{\Gamma }$ contains an
element of this $\omega _{\Gamma }$--convex neighborhood system.  Then
a Moser type argument shows that any small enough neighborhood $\nu
C_\Gamma \subset X_\Gamma$ is symplectomorphic to a neighborhood $\nu
C$ of $C$ in $(X,\omega)$, and hence $\nu C$ contains an
$\omega$--convex neighborhood $U_C$.  In the construction of
$(X_\Gamma,\omega_\Gamma)$ we will use simple models for the surfaces
which are symbolized by the vertices of the plumbing graph $\Gamma$
(similarly to the approach we applied for the central vertex of a
starshaped graph in \cite{GaS}) and will apply a toric construction
for the edges of $\Gamma$ (similarly to the construction along the
legs in \cite{GaS}).  Since this construction might be of independent
interest, we state it as

\begin{thm}\label{t:convneighb}
If $C=C_1\cup \ldots \cup C_n\subset (X, \omega )$ is a collection of
symplectic surfaces in a symplectic 4--manifold $(X, \omega )$
intersecting each other $\omega$--orthogonally according to the
negative definite plumbing graph $\Gamma$ and $\nu C \subset X$ is an open
set containing $C$, then $C$ admits an $\omega$--convex neighborhood
$U_C\subset \nu C \subset (X, \omega )$. In particular, the complement
$X-{\rm {int }}U_C$ is a strong concave filling of its contact boundary.
\end{thm}

\begin{rem}
{\rm Using Grauert's result \cite{grau} it is not hard to show that
  $C$ admits a neighborhood which is a weak symplectic filling of an
  appropriate contact structure on its boundary. (A weak filling is
  one where the symplectic structure is positive on the contact planes
  on the boundary, as opposed to a strong filling, where the contact
  structure is induced by a Liouville vector field transverse to the
  boundary.) Therefore the complement of this neighborhood is a weak
  concave filling, and although in some cases weak convex fillings can
  be deformed to be strong~\cite{Elistrong}, no similar result for
  concave fillings is known. Weak fillings, however, are not suitable
  for the gluing constructions we will apply later, since in the weak
  case the contact structures do not determine the behavior of the
  symplectic forms near the boundaries. In the strong case, the
  Liouville vector fields allow us to glue symplectic forms when the
  contact forms agree.  Hence we verify the existence of an
  $\omega$--convex neighborhood, providing the desired strong concave
  filling of the boundary of the appropriate neighborhood.  Notice
  also that in this first step the further assumptions on the plumbing
  graph $\Gamma$ (listed in (1) and (2) of Theorem~\ref{t:main}) are
  not necessary.}
\end{rem}

After finding the $\omega$--convex neighborhood $U_C \subset (X,
\omega ) $ we would like to compare the induced contact structure $\xi
_C$ on $\partial U_C$ to the Milnor fillable contact structure $\xi
_{M}$ on $\partial Z_1$ (given as the 2--plane field of complex
tangencies on the link).  To this end we describe an open book
decomposition of $\xi_C$ and (using a result of \cite{CNP}) relate it
to an open book decomposition of the Milnor fillable contact structure
$\xi _{M}$. A natural open book decomposition compatible with $\xi _C$
will be given only under the additional hypothesis that $-s_v-d_v\geq
0$ for each vertex $v$ of $\Gamma$, and the relation to some open book
decomposition compatible with $\xi _M$ will be established in the
two cases listed by Theorem~\ref{t:main}. It is natural to conjecture,
however, that these further technical assumptions are unnecessary, hence
we state

\begin{conj}\label{c:main}
The contact structures $\xi _C$ and $\xi _M$ are contactomorphic for
any negative definite plumbing graph $\Gamma$, consequently the
symplectic structure $\omega _Z$ on the 4--manifold $Z$ of
Theorem~\ref{t:main} exists for any negative plumbing graph $\Gamma$.
\end{conj}

The paper is organized as follows: In Section~\ref{sec:gen} we recall
some basics of normal surface singularities.  Section~\ref{sec:second}
is devoted to the description of the $\omega$--convex neighborhoods of
the configuration $C\subset (X, \omega )$ and hence the proof of
Theorem~\ref{t:convneighb}.  In Section~\ref{sec:third}, under the
additional assumption $-s_v-d_v\geq 0$ mentioned above, we describe an
open book decomposition of $(U_C, \xi _C) $ compatible with the
contact structure induced on the boundary of the $\omega$--convex
neighborhood, while in Section~\ref{sec:fourth} we prove
Theorem~\ref{t:main}.

{\bf Acknowledgements}: The second author was partially supported by EU
Marie Curie TOK project BudAlgGeo and by OTKA T49449. Both authors
wish to acknowledge support by ZA-15/2006 Bilateral Project (South
African NRF Grant number 62124).  The second author also would like to
thank Andr\'as N\'emethi and S\'andor Kov\'acs for helpful
discussions.

\section{Generalities on normal surface singularities} \label{S:generalities}
\label{sec:gen}
For the sake of completeness, in this section we collect 
some of the basic results regarding normal surface singularities.
For general reference see \cite{lauf,lw, nem, wahl}.

A \emph{complex germ} $(V,0)$ is an equivalence class of subsets of
$\bfc^n$, where two subsets are equivalent if they agree on some open
neighborhood of $0$. A \emph{germ} $f \colon (\bfc^n,0) \to (\bfc,0)$
of a holomorphic function is an equivalence class of holomorphic
functions from $(\bfc^n,0)$ to $(\bfc,0)$, where two functions are
equivalent if they agree on some open neighborhood of $0 \in
\bfc^n$. Note that the ``inverse image of $0$'' under a germ of a
holomorphic function is naturally a complex germ. Also note that all
derivatives of a holomorphic germ are well defined at $0$.  The
complex germ $(V, 0)$ is a \emph{surface  singularity} if there are
germs of holomorphic functions $f_i\colon (\bfc ^n , 0)\to (\bfc , 0)$ 
($i=1, \ldots , m$) such that
\begin{equation}\label{eq:sing}
(V, 0)=\{ x\in \bfc ^n \mid f_i (x)=0\ \ i=1,\ldots , m \} ,
\end{equation}
and
the rank 
$r(x)$  of the matrix 
\[
( \frac{\partial f_i}{\partial z_j}(x) ) _{i=1, \ldots, m; j=1, \ldots , n}
\] 
is equal to $n-2$ for generic points $x$ of $V$. If $r(x)=n-2$ for all
$x\in V-0$ and $r(0) < n-2$ then the singularity is called
\emph{isolated}.  $(V, 0)$ is \emph{normal} if any bounded holomorphic
function $f \colon V- \{ 0\} \to \bfc $ extends to a holomorphic
function on $V$.  A normal surface singularity is necessarily
isolated. The singularity $(V, 0)$ is a \emph{complete intersection}
if $m=n-2$ in \eqref{eq:sing}, and it is a \emph{hypersurface singularity}
if $n=3$ and $m=1$.

The \emph{link} $L$ of the normal surface singularity $(V,0)$ is defined as the
intersection of $V$ and a sphere $S^{2n-1}_{\epsilon}=\{ x\in \bfc ^n
\mid \vert x\vert =\epsilon\}$. The 3--manifold $L$ is independent of
the embedding of $V$ into $\bfc ^n$, and (provided it is small enough)
independent of
$\epsilon$.

A \emph{resolution} of a singularity $(V, 0)$ is a smooth complex
surface ${\tilde {V}}$ together with a proper holomorphic map $\pi
\colon {\tilde {V}}\to V$ such that $\pi$ restricted to $\pi
^{-1}(V-\{ 0\} )$ is an isomorphism, that is, a diffeomorphism which
is holomorphic in both directions.  The resolution is \emph{good} if
$\pi ^{-1}(0)$ is a normal crossing divisor, that is, in a
decomposition of $\pi ^{-1}(0)=E =E_1\cup \ldots \cup E_k$ into
irreducible components all curves are smooth, intersect each other
transversely and there is no triple intersection.  Such a resolution
always exists, but it is not unique. A resolution is called
\emph{minimal} if it does not contain any rational curve with
self--intersection $(-1)$. The minimal resolution is unique, but might
not be good (in the above sense).  The resolution can be assumed to be
K\"ahler, in such a way that $\pi$ is a biholomorphism away from $0\in
V$. A good resolution can be described by its \emph{dual graph},
where each irreducible component of $E$ is symbolized by a vertex,
each vertex is decorated by the genus and the self--intersection of
the corresponding component, and two vertices are connected if the
corresponding curves intersect each other. Notice that since the
curves $E_i$ are assumed to be smooth, the resulting graph contains no
edge with coinciding endpoints.  It is easy to see that the plumbing
3--manifold defined by the dual graph of a resolution is diffeomorphic
to the link of the singularity at hand.

A resolution graph of a normal surface singularity is always negative
definite, and according to a deep theorem of Grauert \cite{grau}, any
negative definite plumbing graph appears as the graph of a resolution
of an appropriate (and not necessarily unique) normal surface
singularity.  Notice that the link $L$ of the singularity $(V, 0)$
admits a contact structure by considering the complex tangents along
$L$. According to \cite{CNP} this contact structure is unique up to
contactomorphism. It is called the \emph{Milnor fillable} contact
structure on $L$. By a famous result of Bogomolov \cite{Bog} the
complex structure on a resolution ${\tilde {V}}$ can be deformed to a
(possible blow--up of a) Stein filling, hence Milnor fillable contact
structures are necessarily Stein fillable.

A \emph{smoothing} of $(V,0)$ consists of a germ of a complex 3--fold
$({\mathfrak{V}}, 0) $ together with a (germ of a) proper flat
analytic map $f \colon ({\mathfrak {V}}, 0) \to (\Delta , 0)$ (where
$(\Delta , 0) $ is the germ of an open disk in $\bfc$) and an
isomorphism $i\colon (f^{-1}(0), 0) \to (V, 0)$ such that ${\mathfrak
  {V}}-\{ 0 \}$ is nonsingular and $f\vert _{{\mathfrak {V}}-\{ 0\}}$
is a submersion.  By the Ehresman fibration theorem it follows then
that over $\Delta -\{ 0\}$ the map $f$ is a fiber bundle whose fibers
are smooth 2--dimensional Stein manifolds. The typical (nonsingular)
fiber is called the \emph{Milnor fiber} of the smoothing.  Notice that
its boundary is equal to the link of the singularity, and the contact
structure induced on it by the complex tangencies is isotopic to the
Milnor fillable contact structure of the link. Such smoothing does not
necessarily exist for a given singularity; if it does, the Milnor
fiber provides a further Stein filling of the Milnor fillable contact
structure of the link of the singularity.

\section{Construction of $\omega$--convex neighborhoods}
\label{sec:second}

The aim of this section is to prove Theorem~\ref{t:convneighb}.  We
will always assume that $\Gamma$ does not admit an edge from a vertex
back to itself; in other words, the symplectic surfaces $C_i \subset
(X, \omega )$ are assumed to be embedded. The general case involving
immersed surfaces can always be reduced to this situation by
blow--ups.

By applying the following result (which is an application of Moser's
method), the construction of the appropriate neighborhood relies on
constructing model symplectic structures on the plumbing 4--manifold
$X_{\Gamma}$ determined by $\Gamma$. We start with recalling the
Moser-type result.

\begin{thm}[Moser, cf. also \cite{GaS, Sym2}]  \label{t:Moser}
Suppose that $\omega _1$ and $\omega _2$ are symplectic forms on a
$4$--manifold $M$ containing a configuration of smooth surfaces $C =
C_1 \cup \ldots \cup C_n$ which are both $\omega_1$-- and
$\omega_2$--symplectic, with intersections which are both $\omega_1$--
and $\omega_2$--orthogonal. Then $C$ admits symplectomorphic
neighborhoods $(U_1, \omega _1)$ and $(U_2, \omega _2)$ (via a
symplectomorphism which is the identity on $C$) if and only if $\int
_{C_i}\omega _1=\int _{C_i}\omega _2$ for all $i=1, \ldots , n$. \qed
\end{thm}

The rest of the section is occupied by the construction of the model
neighborhoods.  Let $\Gamma$ be a finite graph with vertex set $\{ 1,
2, \ldots, n\}$, with each vertex $v$ labelled with a
self-intersection $s_v \in \Z$, an area $a_v \in \R^+$ and a genus
$g_v \geq 0$.  (As always, $\R^+$ denotes $(0,\infty)$.)  Let $\oma =
(a_1, \ldots, a_n)^T \in (\R^+)^n$. Assume that $\Gamma$ has no edges
from a vertex back to itself. Let $Q$ be the associated $n \times n$
intersection matrix for $\Gamma$, so that $Q_{ii} = s_i$ and $Q_{ij}$
is the number of edges from vertex $i$ to vertex $j$.  (Notice that
the off--diagonals of $Q$ are therefore all nonnegative.) The result
we will prove will be slightly more general than needed because we
will assume a condition more general than that $Q$ is negative
definite.

In~\cite{GaS} we defined a {\em neighborhood $5$-tuple} as a $5$-tuple
$(X,\omega,C,f,V)$ such that $(X,\omega)$ is a symplectic
$4$-manifold, $C$ is a collection of symplectic surfaces in $X$
intersecting $\omega$-orthogonally, $f\colon X \to [0,\infty)$ is a
smooth function with no critical values in $(0,\infty)$ and with
$f^{-1}(0) = C$, and $V$ is a Liouville vector field on $X - C$ with
$df(V) > 0$. From this it easily follows that, for small $t>0$,
$f^{-1}[0,t]$ is an $\omega$--convex tubular neighborhood of $C$.

\begin{prop} \label{P:nbhd5tuple} If there exists a vector $\omz \in (\R^+)^n$
  with $-Q \omz = \frac{1}{2\pi} \oma$ then there exists a
  neighborhood $5$-tuple $(X,\omega,f,C,V)$ such that $C$ is a
  configuration of symplectic surfaces $C_1 \cup \ldots \cup C_n$
  intersecting $\omega$-orthogonally according to the graph $\Gamma$,
  with $C_i \cdot C_i = s_i$, $\int_{C_i} \omega = a_i$ and
  $\genus(C_i) = g_i$.
\end{prop}

Before giving the proof we give a quick survey of the necessary facts
about toric moment maps on symplectic $4$-manifolds. These results are
all standard except that here we suppress the importance of the torus
action and focus instead on how the geometry of the moment map image
determines the smooth and symplectic topology of the total space; from
a $4$-manifold topologist's point of view a useful exposition can be
found in~\cite{Sym3}. Suppose that $\mu \colon X \to \R^2$ is a toric
moment map on a symplectic $4$--manifold $(X,\omega)$ with connected
fibers and with $\partial X = \emptyset$.
\begin{enumerate}
\item Associated to $\mu$ we have coordinates $(p_1,q_1,p_2,q_2)$ on
$X$, with $p_i \in \R$ and $q_i \in \R/2\pi\Z$, such that
$\mu(p_1,q_1,p_2,q_2)=(p_1,p_2)$ and $\omega = dp_1 \wedge dq_1 + dp_2
\wedge dq_2$.
\item The image $\mu(X) \subset \R^2$ has
polygonal boundary with edges of rational slope. Where two edges with
primitive integral tangent vectors $(a,b)^T$ and $(c,d)^T$ (oriented
by $\partial \mu(X)$) meet at a vertex, we have the ``Delzant
condition'':
\[ \det \left( \begin{array}{cc} a & c \\ b & d \end{array} \right) =
1. \]
\item The fibers over interior points of $\mu(X)$ are tori (with coordinates
  $(q_1,q_2)$). The fiber above a point in the interior of an edge of
  $\partial \mu(X)$ with primitive integral tangent vector $(a,b)^T$ is a
  circle with coordinate $a q_1 + b q_2$, so that the $(-b,a)$-circles in a
  nearby $(q_1,q_2)$-torus bound disks. The fiber above a vertex of $\partial
  \mu(X)$ is a single point.
\item Any other symplectic $4$--manifold $(X',\omega')$ with toric
moment map $\mu'\colon X' \to \R^2$ with connected fibers and with
$\mu'(X') = \mu(X)$ is symplectomorphic to $(X,\omega)$ via a
fiber-preserving symplectomorphism.                                                                           Furthermore, the closure of any $2$-dimensional submanifold $B$ of $\R^2$ that has a rational slope
polygonal boundary satisfying the Delzant conditions occurs as
the image of a toric moment map on some symplectic $4$-manifold (with
connected fibers).
\item \label{I:GL2Z} 
Given any matrix $A \in GL(2,\Z)$, there exists a toric moment map
$\mu_A\colon (X,\omega) \to \R^2$ such that $\mu_A(X) = A\mu(X)$ and
such that the coordinates $(p_1',q_1',p_2',q_2')$ associated to
$\mu_A$ are related to the coordinates $(p_1,q_1,p_2,q_2)$ associated
to $\mu$ via the following transformation:
\[ \left( \begin{array}{c} p_1' \\ p_2' \end{array} \right) = A \left(
\begin{array}{c} p_1 \\ p_2 \end{array} \right) ,  \;  
\left( \begin{array}{c} q_1' \\ q_2' \end{array} \right) = A^{-T} \left(
\begin{array}{c} q_1 \\ q_2 \end{array} \right) .
\]
(Here $A^{-T} = (A^{-1})^T$.)
\item  \label{I:Liouville} The vector field $x \partial_x + y \partial_y$ radiating out
from the origin in $\R^2$ lifts to a Liouville vector field $V = p_1
\partial_{p_1} + p_2 \partial_{p_2}$ on $X - \mu^{-1}(\partial
\mu(X))$. Given some $A \in GL(2,\Z)$, the change of coordinates
discussed in the preceding point transforms $V$ to $V' = p'_1
\partial_{p'_1} + p'_2 \partial_{p'_2}$.
\item \label{I:annulusend} Looking at a very specific case, if $R = (x_0,x_1)
  \times [y_0,y_1)$ is an open subset of $B=\mu(X)$ (hence
  $(x_0,x_1)\times \{ y_0\}\subset \partial B$), then the set
  $\mu^{-1}(R)$ is diffeomorphic to $(x_0,x_1) \times S^1 \times
  D^2_\rho$, where $D^2_\rho$ is an open disk in $\R^2$ of radius
  $\rho = \sqrt{2(y_1-y_0)}$ centered at the origin.  Furthermore,
  $\omega|_{\mu^{-1}(R)} = dt \wedge d\alpha + r dr \wedge d\theta$, where $t
  \in (x_0,x_1)$, $\alpha \in \R / 2 \pi \Z$ and $(r,\theta)$ are standard polar
  coordinates on $D^2_\rho$, and with these coordinates,
  $\mu(t,\alpha,r,\theta) = (t,\frac{1}{2} r^2 + y_0)$, i.e. $p_1 =
  t$, $q_1 = \alpha$, $p_2 = \frac{1}{2} r^2 + y_0$, $q_2 =
  \theta$. Then $\mu^{-1}(\partial R) = \mu^{-1}((x_0,x_1) \times \{
  y_0\} )$ is a cylinder $(x_0,x_1) \times S^1 \times \{0\}$ with
  symplectic area $2\pi(x_1-x_0)$. The Liouville vector field $p_1
  \partial_{p_1} + p_2
\partial_{p_2}$ then becomes $V = t \partial_t + (\frac{1}{2} r +
\frac{y_0}{r}) \partial_r$. (Note that $V$ is clearly undefined at
$r=0$ except in the special case that $y_0 = 0$.)
\end{enumerate}

\begin{proof}[Proof of Proposition~\ref{P:nbhd5tuple}]
  Fix a vector $\omz = (z_1, \ldots, z_n)^T \in (\R^+)^n$ with $-Q \omz =
  \frac{1}{2\pi} \oma$. For each vertex $v$ and for each edge $e$ meeting
  $v$, choose an integer $s_{v,e}$ such that $\sum s_{v,e} = s_v$, where this
  sum and other similar sums below are taken over all edges meeting the given
  vertex $v$. Also, for each vertex $v$ and each edge $e$ meeting $v$, letting
  $w$ be the vertex at the other end of $e$, let $x_{v,e} = -s_{v,e} z_v -
  z_w$. Note that, for each $v$ we have $\sum x_{v,e} = (-Q \omz)_v =
  \frac{1}{2\pi} a_v > 0$.  Choose a small positive constant $\epsilon$,
  small enough so that for each $v$ we have $\sum (x_{v,e} - \epsilon) > 0$.
  Also choose small positive constants $\delta$ and $\gamma$ satisfying a
  constraint to be stated shortly.

Consider the first quadrant $P = [0,\infty)^2 \subset \R^2$ and let $g \colon
P \ra [0,\infty)$ be a smooth function satisfying the following
properties (see Figure~\ref{F:g}):
\begin{enumerate}
\item $0$ is the only critical value of $g$.
\item $g^{-1}(0) = \partial P$.
\item If $y-x \geq \gamma$ then $g(x,y) = x$.
\item If $y-x \leq - \gamma$ then $g(x,y) = y$.
\item For all $x,y$ we have $g(x,y) = g(y,x)$.
\item In the region $-\gamma \leq y-x \leq \gamma$, the level sets
$g^{-1}(t)$, for $t > 0$, are smooth curves symmetric about the line
$y=x$, with slope changing monotonically as a function of $y-x$ from
$0$ to $\infty$.
\end{enumerate} 
\begin{figure}[ht!]
\labellist
\small\hair 2pt
\pinlabel $y$ at 3 324
\pinlabel $x$ at 324 3
\pinlabel $\gamma$ at -12 60
\pinlabel $\gamma$ at 60 -12
\pinlabel {$y-x > \gamma$} at 104 255
\pinlabel {$y-x < -\gamma$} at 255 104
\endlabellist
\centering
\includegraphics[scale=0.4]{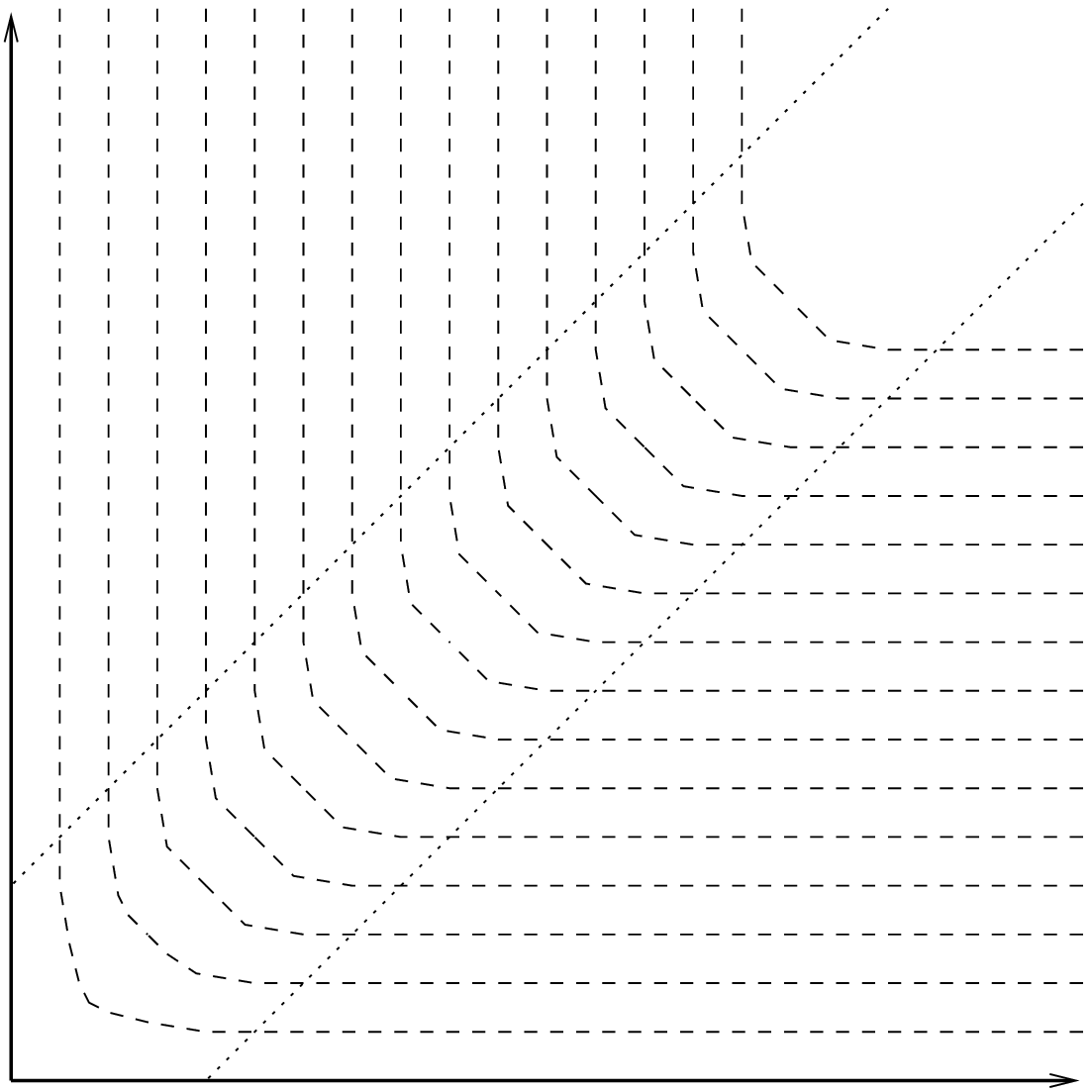}
\caption{Contour plot of $g$.}
\label{F:g}
\end{figure}

The constants $\delta$ and $\gamma$ should satisfy the
following constraint: For each vertex $v$ and for each edge $e$ incident
to $v$, the line passing through $(0,\epsilon)$ with tangent vector
$(1,-s_{v,e})$ should intersect $g^{-1}(\delta)$ in the region $y-x >
\gamma$. By symmetry we will also have that the line passing through
$(\epsilon,0)$ with tangent vector $(-s_{v,e},1)$ intersects
$g^{-1}(\delta)$ in the region $y-x < -\gamma$. Note that if $s_{v,e} <0$, this
constraint is simply the constraint that $\gamma < \epsilon$.

For each edge $e$ we now construct a neighborhood $5$-tuple $(X_e,
\omega_e, f_e, C_e, V_e)$ as follows (see Figure~\ref{F:ge}): Consider
the two vertices at the ends of $e$ and arbitrarily label one $v$ and
the other $v'$. Let $g_e(x,y) = g(x - z_v, y - z_{v'})$, a function
from $P+(z_v,z_{v'})$ to $[0,\infty)$. Let $R_e$ be the open subset of
$g_e^{-1}[0,\delta)$ between the line passing through
$(z_v,z_{v'}+2\epsilon)$ with tangent vector $(1,-s_{v,e})$ and the
line passing through $(z_v+2\epsilon,z_{v'})$ with tangent vector
$(-s_{v',e},1)$.  Let $(X_e,\omega_e)$ be the unique connected
symplectic $4$-manifold with toric moment map $\mu_e\colon X_e \ra
\R^2$ such that $\mu_e(X_e) = R_e$. Let $C_e = \mu_e^{-1}(\partial
R_e)$, $f_e = g_e \circ \mu_e$ and let $V_e$ be the Liouville vector
field obtained by lifting the radial vector field emanating from the
origin in $\R^2$, as in item~(\ref{I:Liouville}) in the discussion of
toric geometry above. Note that $df_e(V_e) > 0$ because $dg_e(x
\partial_x + y
\partial_y) > 0$, which is true because $z_v > 0$ and $z_{v'} > 0$.
(Topologically,
$C_e$ is just a union of two disks meeting transversely at one
point and  $X_e$ is a $4$--ball neighborhood of $C_e$.)
\begin{figure}[ht!]
\labellist
\small\hair 2pt
\pinlabel {$(z_v,z_{v'})$} at 75 130
\pinlabel {$(z_v,z_{v'}+\epsilon)$} [r] at 85 285
\pinlabel {$(z_v,z_{v'}+2\epsilon)$} [r] at 85 426
\pinlabel {$(z_v+\epsilon,z_{v'})$} at 223 126
\pinlabel {$(z_v+2\epsilon,z_{v'})$} at 368 126
\pinlabel {$R_{e,v}$} at 109 353
\pinlabel {$R_{e,{v'}}$} at 311 165
\pinlabel $y$ [b] at 29 454
\pinlabel $x$ [l] at 397 30
\pinlabel {$g_e^{-1}(\delta)$} [l] at 447 186
\endlabellist
\centering
\includegraphics[scale=0.5]{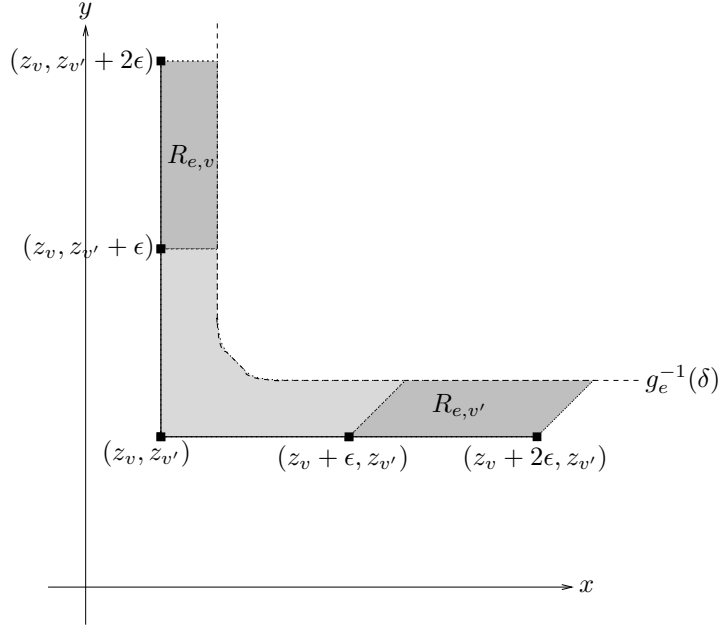}
\caption{The moment map image $R_e$ of $(X_e,\omega_e)$; in this
example $s_{v,e} = 0$ and $s_{v',e} = -1$.}
\label{F:ge}
\end{figure}

Also let $R_{e,v}$ be the open subset of $R_e$ between the parallel
lines passing through $(z_v,z_{v'}+\epsilon)$ and
$(z_v,z_{v'}+2\epsilon)$ with tangent vector $(1,-s_{v,e})$, and let
$R_{e,v'}$ be the open subset of $R_e$ between the parallel lines
passing through $(z_v+\epsilon,z_{v'})$ and $(z_v+2\epsilon,z_{v'})$
with tangent vector$(-s_{v',e},1)$. By the constraints on $\delta$ and
$\gamma$, these are both parallelograms, open on three sides.

Now we introduce two reparametrizations of this neighborhood
$5$-tuple, one for each of the vertices $v$ and $v'$, using matrices
$A_v, A_{v'} \in GL(2,\Z)$ as in item~(\ref{I:GL2Z}) preceding this
proof. These matrices are:
\[ A_v = \left( \begin{array}{cc} -s_{v,e} & -1 \\ 1 & 0 \end{array}
\right), \; A_{v'} = \left( \begin{array}{cc} -1 & -s_{v',e} \\ 0 & 1
\end{array} \right)  . \] 
The reader should at this point verify that $A_v$ transforms $R_{e,v}$
into the region $(x_{v,e} - 2\epsilon, x_{v,e} -
\epsilon) \times [z_v, z_v+\delta)$ and that $A_{v'}$ transforms
$R_{e,v'}$ into the region $(x_{v',e} - 2\epsilon, x_{v',e} -
\epsilon) \times [z_{v'}, z_{v'}+\delta)$. Referring to
item~(\ref{I:annulusend}) in the toric discussion preceding this proof,
we see that on $\mu_e^{-1}(R_{e,v})$ and on $\mu_e^{-1}(R_{e,v'})$ we
can write everything down in particularly nice local coordinates as
follows: On $\mu_e^{-1}(R_{e,v})$ we have:
\begin{itemize}
\item $\mu_e^{-1}(R_{e,v}) \cong (x_{v,e}-2\epsilon,x_{v,e} - \epsilon)
\times S^1 \times D^2_{\sqrt{2\delta}}$ with corresponding coordinates
$(t,\alpha,r,\theta)$.
\item In these coordinates, $\omega_e = dt \wedge d\alpha + r dr
\wedge d\theta$.
\item $C_e \cap \mu_e^{-1}(R_{e,v}) = (x_{v,e}-2\epsilon,x_{v,e} - \epsilon)
\times S^1 \times \{0\}$.
\item $f_e = \frac{1}{2} r^2$. 
\item $V_e = t \partial_t + (\frac{1}{2} r + \frac{z_v}{r} ) \partial_r$.
\end{itemize}
On $\mu_e^{-1}(R_{e,v'})$ we have exactly the same formulae but with
each occurrence of $v$ replaced with $v'$. 

Now we will construct neighborhood $5$-tuples associated to the
vertices so that they can be glued to the neighborhoods constructed
above using the explicit coordinates that we have just seen in the
preceding paragraph.  Lemma~2.4 from~\cite{GaS} tells us that for each
vertex $v$ we can find a compact surface $\Sigma_v$ of genus $g_v$
with a symplectic form $\beta_v$ and Liouville vector field $W_v$
($\beta_v$ and $W_v$ both defined on all of $\Sigma_v$) such that
$\Sigma_v$ has one boundary component $\partial_{e,v} \Sigma_v$ for
each edge $e$ incident with $v$ and such that there exists a collar
neighborhood $N_{e,v}$ of each $\partial_{e,v} \Sigma_v$ parametrized
as $(x_{v,e}-2\epsilon,x_{v,e}-\epsilon] \times S^1$ on which $\beta_v
= dt \wedge d\alpha$ and $W_v = t\partial_t$. (Here we use the
constraint we imposed on $\epsilon$, namely that, for each vertex $v$
we have $\sum (x_{v,e} - \epsilon) > 0$.)  Note that $\int_{\Sigma_v}
\beta_v = 2 \pi \sum (x_{v,e} - \epsilon)$. Then our neighborhood
$5$-tuple for the vertex $v$ is:
\begin{eqnarray*} (X_v &=& (\Sigma_v - \partial \Sigma_v)
\times D^2_{\sqrt{2\delta}}, \\
\omega_v &=& \beta_v + r dr \wedge d\theta, \\
C_v &=& \Sigma_v - \partial \Sigma_v, \\
f_v &=& \frac{1}{2} r^2, \\
V_v &=& W_v + (\frac{1}{2} r + \frac{z_v}{r}) \partial_r). 
\end{eqnarray*}

These neighborhoods can then be glued to the neighborhoods for the edges as
follows: For each edge $e$ with incident vertices $v$ and $v'$, glue the end
$(N_{e,v} - \partial_{e,v} \Sigma_v) \times D^2_{\sqrt{2 \delta}}$ of
$X_v$ to the end $\mu_e^{-1}(R_{e,v})$ of $X_e$ by identifying the
$(t,\alpha,r,\theta)$ coordinates, and similarly glue $(N_{e,v'} -
\partial_{e,v'} \Sigma_v') \times D^2_{\sqrt{2 \delta}}$ to
$\mu_e^{-1}(R_{e,v'})$. The result is the $5$-tuple
$(X,\omega,C,f,V)$. 

We now verify that the areas and self-intersections of the surfaces in
$C$ are correct.  For the areas, note that the closed surface $C_v
\subset X$ is
the union of $(\Sigma_v - \partial \Sigma_v) \times 0$ in $X_v$ with
the various disks $\mu_e^{-1}(\partial_v R_e) \subset X_e$, where $\partial_v R_e$
is one of the two edges making up $\partial R_e$. The area of
$(\Sigma_v - \partial \Sigma_v) \times 0$ is $2\pi \sum (x_{v,e} -
\epsilon)$, the area of each disk is $2\pi (2\epsilon)$ and the area
of each overlapping cylinder is $2\pi \epsilon$, so the total area is
$2\pi \sum x_{v,e} = a_v$. For the self--intersections, note that the
boundary of a tubular neighborhood of $C_v$ is a $3$-manifold homeomorphic to
$\Sigma_v \times S^1$ with the boundary components Dehn filled with solid
tori. Looking at how the matrices $A_v$ (or $A_{v'}$) transform the regions
$R_e$, and following the argument at the end of the proof of
\cite[Proposition~2.3]{GaS}, we see that the $(1,s_{v,e})$ curves in each
$\partial_{v,e} \Sigma_v \times S^1$ are filled in by disks. So this
$3$--manifold is the $S^1$--bundle over $C_v$ of Euler class $\sum s_{v,e} =
s_v$.
\end{proof}

In order to apply Proposition~\ref{P:nbhd5tuple} in the proof of
Theorem~\ref{t:Moser} we need to show that the symmetric matrix $Q$ defined by
the graph $\Gamma$ of the symplectic surfaces $C_1 \cup \ldots \cup C_n
\subset (X, \omega )$ satisfies the property that the equation
\[
-Q\omz = \frac{1}{2\pi }\oma
\]
admits a solution $\omz = (z_1, \ldots , z_n) \in (\bfr ^+)^n$ for any given
$\oma \in (\bfr ^+)^n$. The basis of our argument is the following simple
linear algebra observation:
\begin{lem}\label{l:linalg}
  Suppose that the bilinear form $(x,y)$ is given by the negative definite
  symmetric matrix $Q$ with only nonnegative off--diagonals in the basis $\{
  E_i\}$. If for a vector $x$ the inequalities $(x,E_i)\leq 0$  ($i=1,
  \ldots , n$) are all satisfied, then all coordinates of $x$ are nonnegative.
\end{lem}
\begin{proof}
  Let us expand $x$ in the basis $\{ E_i\}$ and denote the resulting
  $n$--tuple by $x$ as well.  Suppose that $x=x_1-x_2$ where $x_i$ has only
  nonnegative entries for $i=1,2$, and the supports of $x_1$ and $x_2$ are
  disjoint. Take $E_i$ from the support of $x_2$.  Then by the assumption
\[
(x,E_i)=(x_1,E_i)-(x_2, E_i)\leq 0
\]
implying that $(x_1, E_i)\leq (x_2, E_i)$. Summing for all basis vectors $E_i$
in the support of $x_2$ and multipling the inequalities with the positive
coefficients they have in $x_2$ we get
\[
(x_1,x_2) \leq (x_2, x_2).
\]
Since the supports of $x_1$ and $x_2$ is disjoint (and the
off--diagonals in $Q$ are all nonnegative, that is, $(E_i, E_j) \geq
0$ once $i\neq j$), we have that $(x_1, x_2)\geq 0$. On the other hand, $Q$ is
negative definite, so $(x_2, x_2)\leq 0$.  This implies that $(x_2,
x_2)=0$, which by definiteness implies that $x_2=0$, hence $x=x_1$, verifying 
the lemma.
\end{proof}

\begin{cor}
For any $\oma \in (\bfr ^+)^n$ the vector $-Q^{-1}\oma $ is in $(\bfr ^+)^n$.
\end{cor}
\begin{proof}
Suppose that $\oma $ is in $ (\bfr ^+)^n$ and consider 
$\omb =-Q^{-1}\oma $. Then $-\oma = Q\omb$ is a vector with
only nonpositive coordinates, that is, $(\omb , E_i )\leq 0$ for all $i$.
The application of Lemma~\ref{l:linalg} then finishes the proof.
\end{proof}

\begin{proof}[Proof of Theorem~\ref{t:convneighb}]
  By the above corollary and Proposition~\ref{P:nbhd5tuple}, there
  exists a neighborhood $5$-tuple $(X_\Gamma, \omega_\Gamma, f_\Gamma,
  C_\Gamma, V_\Gamma)$ for the given plumbing graph $\Gamma$
  (decorated with $a_i =\int _{C_i} \omega$). By basic results in
  differential topology, there exists an open neighborhood $U$ of $C$
  in $X$ which is diffeomorphic to $f_\Gamma^{-1}(t)$ for some small
  $t > 0$, via a diffeomorphism sending $C$ to $C_\Gamma$. By
  Theorem~\ref{t:Moser}, we can make this diffeomorphism into a
  symplectomorphism, after possibly taking a smaller neighborhood of
  $C$ and a smaller value for $t$. Since in the neighborhood 5--tuple
  every neighborhood of $C_{\Gamma }$ contains an $\omega
  _{\Gamma}$--convex neighborhood, its image under the
  symplectomorphism provides $U_C\subset (X, \omega )$.
\end{proof}

\section{Open book decompositions on $\partial U_C$}
\label{sec:third}
Suppose that the plumbing graph $\Gamma$ satisfies the additional hypothesis
that for every vertex $v$ the self--intersection (homological square) $s_v$
and the valency $d_v$ 
\[
-s_v-d_v\geq 0
\]
holds. In this section we describe an open book decomposition on
$\partial U_C$ compatible with the contact structure induced on it as
an $\omega$--convex neighborhood of $C$.  We begin with a lemma about
``open book decompositions'' (OBDs) on $3$--manifolds with
boundary. By an OBD on a $3$--manifold $M$ with $\partial M \neq
\emptyset$ we mean a pair $(B,\pi)$, where $B \subset M - \partial M$
is a link and $\pi \colon M - B \to S^1$ is a fibration which behaves
as open books usually behave near $B$ and which restricts to $\partial
M$ to give an honest fibration of $\partial M$ over $S^1$. When the
pages are oriented, this induces an orientation on $B$ as the boundary
of a page.

\begin{lem} \label{L:OBDbuildingblock} Consider $M = [0,1] \times S^1 \times
  S^1$ with coordinates $t \in [0,1]$, $\alpha \in S^1$ and $\beta \in
  S^1$.  Given a nonnegative integer $m$ there exists an OBD $(B,\pi)$
  on $M$ such that the following conditions hold:
\begin{enumerate}
\item $\pi|_{\{ 0\} \times S^1 \times S^1} =
  \beta$
\item $\pi|_{\{ 1\} \times S^1 \times S^1} = \beta + m\alpha$
\item The pages $\pi^{-1}(\theta)$
  are transverse to $\partial_\beta$.
\item The binding $B$ is tangent to $\partial_\beta$. 
\item $B$ has $m$ components $B_1,
  \ldots, B_{m}$, which we can take to be $B_i = \{1/2\} \times
  \{(2\pi i)/m\} \times S^1$. 
\item When the pages are oriented so that
  $\partial_\beta$ is positively transverse then $B_1, \ldots, B_m$
  are oriented in the positive $\partial_\beta$ direction.
\end{enumerate}
\end{lem}

\begin{proof}
 If $m=0$ we use the map $\pi = \beta$ on all of $M$ and have
 $B=\emptyset$. Otherwise the proof follows directly from the
 following observation which we leave to the reader to verify (with
 the aid of Figure~\ref{F:OBDbuildingblock}): Consider $P =[0,1]
 \times [0,1] \times S^1$ with coordinates $(x,y,\theta)$. There is an
 OBD $(B_P,\pi_P)$ on $P$ with $B_P = \{1/2\} \times \{1/2\} \times
 S^1$, such that $f|_{\{0\} \times [0,1] \times S^1} = \theta$,
 $f|_{[0,1] \times \{0\} \times S^1} = \theta$, $f|_{[0,1] \times
   \{1\} \times S^1} = \theta$ and $f|_{\{1\} \times [0,1] \times S^1}
 = \theta + 2\pi y$. When the pages are oriented so that
 $\partial_\theta$ is positively transverse then $B_P$ is oriented in
 the positive $\partial_\theta$ direction.
\begin{figure}[ht!]
\labellist
\small\hair 2pt
\pinlabel $x$ at 14 31
\pinlabel $y$ at 183 51
\pinlabel $\theta$ at 78 177
\pinlabel $B_P$ at 111 59
\endlabellist
\centering
\includegraphics[scale=0.85]{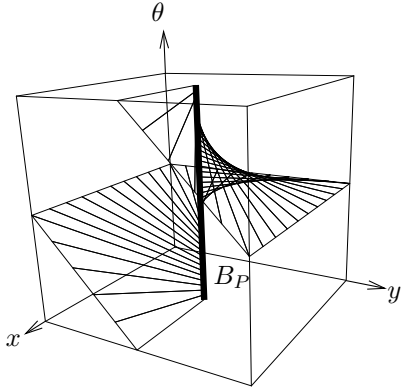}
\caption{Building block for OBD's. The shaded surface indicates a page.}
\label{F:OBDbuildingblock}
\end{figure}
Given this observation, the lemma can be proved by stacking $m$ of the
above models side-by-side (in the $y$ direction). Some trivial smoothing is required, of course.
\end{proof}

Recall that a plumbed $3$-manifold $M=M_{\Gamma }$ constructed
according to a plumbing graph $\Gamma$ decomposes along a collection
of tori $\{T_e\}$, indexed by the edges of $\Gamma$, into
codimension-$0$ pieces $\{M_v\}$, indexed by the vertices of
$\Gamma$. Each $M_v$ fibers overs a compact surface $\Sigma_v$ with
each boundary component $\partial_{v,e} M_v$ of $M_v$ fibering over a
corresponding boundary component $\partial_{v,e} \Sigma_v$ of
$\Sigma_v$. On each torus $T_e$ there are thus two fibrations over
$S^1$, coming from the vertices at the two ends of $e$. We say that an
OBD on $M$ is \emph{horizontal} if the pages are transverse to the
fibers on each $M_v$ and transverse to both types of fibers on each
$T_e$ and if the binding components are disjoint from the $T_e$'s and
are fibers of the fibration of the corresponding $M_v$'s. (Note that
this definition depends on identifying $M$ as a plumbed $3$-manifold
and specifying the fibrations on each $M_v$.) In addition, we can
orient the binding components as boundary components of a page, with
the page oriented so as to intersect fibers positively; we require
this orientation to point in the positive fiber direction.  (For more
about horizontal OBD's, see \cite{TBur}.)

Now we refer to the notation of Proposition~\ref{P:nbhd5tuple} and its proof.
For any small enough $t>0$, $M = f^{-1}(t)$ is a plumbed $3$-manifold. We may
take the separating tori $\{T_e\}$ to be $T_e = \mu_e^{-1}(g_e^{-1}(t) \cap
L)$, where $L$ is the line $(y-z_{v'})-(x-z_v) = 0$ in $R_e$. Let $\xi _C= \ker
(\imath_V \omega|_M)$ be the contact structure induced on $M$ by the Liouville
vector field $V$ and the symplectic structure $\omega$. 

\begin{prop}\label{p:openbook}
Suppose that the plumbing graph $\Gamma$ satisfies the additional
hypothesis that $p_v= -s_v - d_v$ is nonnegative for every vertex $v$
of $\Gamma$.  Then there exists a horizontal OBD on $M$ supporting
$\xi$ with $p_v$ binding components in each fibered piece $M_v$.  This
OBD is independent of the areas $a_1, \ldots, a_n$ of the symplectic
surfaces $C_1, \ldots, C_n$, and therefore the various contact
structures induced by the different symplectic structures for
different $\oma \in (\bfr ^+)^n$ are all isotopic.
\end{prop}

\begin{proof}
Referring to the proof of Proposition~\ref{P:nbhd5tuple}, we see that
$M$ is built by gluing the $f_v^{-1}(t)$'s to the
$f_e^{-1}(t)$'s. Recall that $f_v^{-1}(t) = (\Sigma_v - \partial
\Sigma_v) \times S^1_{\rho}$, where $S^1_{\rho}$ is the circle of
radius $\rho = \sqrt{2t}$. Each $f_e^{-1}(t)$ is a submanifold of
$X_e$ which has toric coordinates $(p_1,q_1,p_2,q_2)$.  The OBD we
construct will be the $S^1_{\rho}$ coordinate function $\theta$ on
each $f_v^{-1}(t)$ and the function $q_1 + q_2$ on each
$f_e^{-1}(t)$. We will put in binding components in the $(x_{v,e} -
2\epsilon,x_{v,e} - \epsilon) \times S^1 \times S^1_{\rho}$ overlaps
where the gluing happens, in order to ``interpolate'' from $\theta$ to
$q_1 + q_2$. In order to do this, we must transform the function
$q_1+q_2$ into the $(t,\alpha,\theta)$ coordinates on each $(x_{v,e} -
2\epsilon,x_{v,e} - \epsilon) \times S^1 \times S^1_{\rho}$ and
$(x_{v',e} - 2\epsilon,x_{v',e} - \epsilon) \times S^1 \times
S^1_{\rho}$, using the transformations given by the matrices $A_v$ and
$A_{v'}$. We see that the change of coordinates associated with $A_v$
at the end $R_{e,v}$, transforms $q_1 + q_2$ into the function
$(-s_{v,e}-1) \alpha + \theta$ and that the change associated with
$A_{v'}$ transforms $q_1 + q_2$ into $(-s_{v',e}-1)\alpha +
\theta$. Thus using Lemma~\ref{L:OBDbuildingblock}, we see that for
each vertex $v$ incident to an edge $e$, if we have nonnegative
integers $p_{v,e}$ with $p_{v,e} = -s_{v,e}-1$ we can interpolate from
$q_1+q_2$ to $\theta$ by introducing $p_{v,e}$ binding components. By
suitably partitioning the $p_v$'s into $p_{v,e}$'s we construct the
desired OBD.

It remains to verify that this OBD is horizontal and supports
$\xi$. The OBD is clearly horizontal on each $f_v^{-1}(t)$ and on the
overlap regions where the binding components are put in. On each
$f_e^{-1}(t)$, we need to see how the fiber directions
$\partial_\theta$ coming from each vertex incident to $e$ transform
via the inverses of the transformations associated to $A_v$ and
$A_{v'}$. This check is straightforward and we see that, at the $v$
end, $\partial_\theta$ becomes $\partial_{q_1}$ and at the $v'$ end,
$\partial_\theta$ becomes $\partial_{q_2}$. Both of these are
transverse to the pages, i.e. the fibers of $q_1 + q_2$.

Lastly, we need to verify that the Reeb vector field for a contact form for $\xi _C$ is
transverse to the pages of this OBD and tangent to the
bindings. However, this is clear because, on $f_v^{-1}(t)$ the Reeb
vector field for the contact form induced by the Liouville vector field is a positive multiple of $\partial_\theta$, and on
$f_e^{-1}(t)$ the Reeb vector field for the contact form induced by the Liouville vector field is a positive multiple of $b_1
\partial_{q_1} + b_2 \partial_{q_2}$ where $dg_e = b_1 dx + b_2 dy$,
and $b_1, b_2 > 0$ by construction of $g_e$.  Notice that in this
construction there was no dependence on the areas $\oma$.
\end{proof}

\section{The proof of Theorem~\ref{t:main}}
\label{sec:fourth}

In order to apply the gluing scheme of symplectic 4--manifolds along
hypersurfaces of contact type (as it is given in \cite{Etn}) we have
to verify that the contact structure $\xi _C$ (given by the toric
picture) and the Milnor fillable contact structure $\xi _{M}$ are
contactomorphic. (Recall that in the previous section we saw that for
a plumbing graph $\Gamma$ for which $-s_v-d_v\geq 0$ holds for every
vertex $v$, the toric approach produces isotopic contact structures
for any input vector $\oma \in (\bfr ^+)^n$.)  In the case of negative
definite starshaped plumbing trees of spheres with three legs this
identification of contact structures relied on the classification of
tight contact structures on certain small Seifert fibered 3--manifolds
\cite{GaS}. Such a classification is not available in
general. Although we strongly believe that the two contact structures
above are contactomorphic in general (which would lead to the
verification of Conjecture~\ref{c:main}), we could prove it only under
strong restrictions on the plumbing graph $\Gamma$, giving the proof
of Theorem~\ref{t:main}.

Recall that each vertex $v$ of the plumbing graph $\Gamma$ is
decorated by two integers: $g_v\geq 0$ denotes the genus of the
surface $\Sigma _v$ corresponding to the vertex $v$, while $s_v$ is
the Euler number of the normal disk bundle of $\Sigma _v$ in the
plumbing 4--manifold $X_{\Gamma }$ (or alternatively the
self--intersection of the homology class $[\Sigma _v]$).  Since
$\Gamma$ is negative definite, we have that $s_v< 0$. As before, let
$d_v$ denote the valency of the vertex $v$, that is, the number of
edges emanating from $v$.  Suppose that $-s_v-d_v\geq 0$ holds for
every vertex $v$.

\begin{proof}[Proof of Theorem~\ref{t:main}]
Suppose first that $\Gamma$ is a tree and for all $v$ we have $g_v=0$,
that is, the surfaces $\Sigma _v$ are all spheres. It is a standard
fact that under these assumptions the boundary 3--manifold of the
plumbing is a rational homology sphere, in which case (according to a
result of Stallings) its binding uniquely determines an OBD.
\cite[Theorem~3.9]{CNP} provides a connection between holomorphic
functions and OBD's compatible with the Milnor fillable contact
structure.  Applying \cite[Theorem~4.1]{CNP} we get an OBD compatible
with $\xi _M$ having the same binding as the OBD we constructed in
Proposition~\ref{p:openbook} (compatible with $\xi _C$). This implies
that $\xi _C$ and $\xi _M$ are contactomorphic in the special case
considered.

If the strict inequality 
\[
-s_v- d_v >2g_v
\]
holds for every vertex, then Proposition~\ref{p:openbook} provides a
horizontal OBD compatible with $\xi _C$ such that it has at least
$2g_i+1$ binding components near every vertex of $\Gamma$.  By
\cite[Theorem~4.1]{CNP} there exists a horizontal OBD compatible with
$\xi _{M}$ which has the same binding as the horizontal OBD
constructed in Proposition~\ref{p:openbook}. Since $2g_i+1>0$ for all
vertices $v_i$ , there are binding components near every vertex. In
this case, however, \cite[Proposition~4.6]{CNP} shows that the two
horizontal OBD's with the same binding are isomorphic, implying that
$\xi _C$ and $\xi _{M}$ are contactomorphic. 

In conclusion, under the assumptions of Theorem~\ref{t:main} the
strong filling $Z_1$ of the Milnor fillable contact link $(Y, \xi _M
)$ and the strong concave filling $X-X_1$ of $(Y, \xi _C)$ have
contactomorphic contact structures on their boundaries, hence the
gluing construction described in \cite{Etn} applies (for a suitably
chosen contactomorphism $\varphi \colon \partial (X-X_1) \to \partial
(-Z_1)$), providing a symplectic stucture on $Z=Z_1\cup _Y (X-X_1)$.
This concludes the proof of the main theorem.
\end{proof}

\end{document}